\documentclass{article}

\usepackage{amsmath}

\title{A FOUR-PARAMETRIC RATIONAL SOLUTION TO PAINLEV\'{E} VI}
\author{Gert Almkvist}

\begin{document}
\maketitle

\textbf{Introduction}

The seminal paper by Okamoto [3] showed how to get a sequence of rational
solutions to Painlev\'{e} VI if you start with a rational seed solution. But
Okamoto did not even write down the B\"{a}cklund transformation. This is
understandable since its denominator is of degree $6$ in $p$ and $q$. Today
we have Maple to handle such things and the author computed hundreds of
examples starting with rational solutions that come from a Riccati equation
and can be expressed by hypergeometric functions (see [2] ). Soon a pattern
emerged. The first $\tau $-function, $\tau _1$, had numerator $1$ and $\tau
_2$ was also rather simple. An explicit formula for $\tau _2$ was found and
proved. As a consequence we have the following main result.

Let $m$ be a positive integer. Define 
\[
W(r,m,s)=\sum_{j=0}^m(-1)^j\binom{r+m+1-j}{m-j}\binom{s+m}jt^{m-j} 
\]
Let further 
\[
T_1=1,T_2=W(r,m,s) 
\]
\[
S_1=1,S_2=W(r-1,m+1,s-1) 
\]
and 
\[
T_{n+1}T_{n-1}=(t^2-t)(T_nT_n^{\prime \prime }-T_n^{\prime \text{ }%
2})+(2t-1)T_nT_n^{\prime }+(n-1)(n+r)T_n^2 
\]
\[
S_{n+1}S_{n-1}=(t^2-t)(S_nS_n^{\prime \prime }-S^{\prime \text{ }%
2})+(2t-1)S_nS_n^{\prime }+(n-1)(n+r-1)S_n^2 
\]
\[
\]
for $n\geq 2$. Then the $T_n$ and $S_n$ are polynomials and 
\[
q_n=t+\frac{t(t-1)}{n+r}\left\{ \frac{S_n^{\prime }}{S_n}-\frac{%
T_{n+1}^{\prime }}{T_{n+1}}-\frac{n+s-1}t-\frac{n+r-m-s}{t-1}\right\} 
\]
solves $P_{VI}$%
\[
q^{\prime \prime }=\frac 12\left( \frac 1q+\frac 1{q-1}\frac 1{q-t}\right)
q^{\prime \text{ }2}-\left( \frac 1t+\frac 1{t-1}+\frac 1{q-t}\right)
q^{\prime }+ 
\]
\[
\frac{q(q-1)(q-t)}{t^2(t-1)^2}\left( \alpha +\frac{\beta t}{q^2}+\frac{%
\gamma (t-1)}{(q-1)^2}+\frac{\delta t(t-1)}{(q-t)^2}\right) 
\]
with 
\[
\alpha =\frac{(n+r)^2}2 
\]
\[
\beta =-\frac{(m+s)^2}2 
\]
\[
\gamma =\frac{(r-s+1)^2}2 
\]
\[
\delta =\frac{1-(n+m)^2}2 
\]
\[
\]

If also $r,s$ are integers then $T_n$ and $S_n$ have integer coefficients,
usually growing very fast with $n.$ One can ask for the smallest integer
coefficients since $q_n$is independent of multiplicative constants in $%
T_{n+1}$ and $S_n.$ As a result we get some intriguing conjectures. E.g.

Let $p\geq 3$ be a prime. Define 
\[
T_1=1,T_2=\frac{W(p,p-1,1)}p 
\]
and $c(n)$ by 
\[
p^2(p+n)(p+n-2)\text{ if }n\equiv 1\text{ mod }p 
\]
\[
\frac{(p+n)(p+n-2)}p\text{ if }n\equiv 0,2\text{ mod }p 
\]
\[
(p+n)(p+n-2)\text{ otherwise} 
\]
Then 
\[
c(n)T_{n+1}T_{n-1}=(t^2-t)(T_nT_n^{\prime \prime }-T_n^{\prime \text{ }%
2})+(2t-1)T_nT_n^{\prime }+(n-1)(n+p)T_n^2 
\]
for $n\geq 2$ has solutions $T_n$ in \textbf{Z}[t] where the coefficients
have no common factor $>1$. We have checked the conjecture for $%
p=3,5,7,11$ up to $n=20$ or more.

In [1] there are other rational solutions which hopefully can be used in the
same way to produce sequences of rational solutions to P$_{VI}.$%
\[
\]
\textbf{1.Deriving the results.}

We follow the notation in Okamoto's paper [3]. Given a solution 
\[
q=q(\mathbf{b}) 
\]
with parameters 
\[
\mathbf{b}=(b_1,b_2,b_3,b_4) 
\]
i.e. $q$ satisfies 
\[
P_{VI}(\frac{(b_3-b_4)^2}2,-\frac{(b_1+b_2)^2}2,\frac{(b_1-b_2)^2}2,-\frac{%
(b_3+b_4)(b_3+b_4+2)}2) 
\]
This equation is equivalent to the system 
\[
\frac{dq}{dt}=\frac{\partial H}{\partial p} 
\]
\[
\frac{dp}{dt}=-\frac{\partial H}{\partial q} 
\]
where 
\begin{eqnarray*}
H &=&H(\mathbf{b})= \\
&&\frac{q(q-1)(q-t)}{t(t-1)}\left\{ p^2-p(\frac{b_1+b_2}q+\frac{b_1-b_2}{q-1}%
+\frac{b_3+b_4}{q-t})+\frac{(b_1+b_3)(b_1+b_4)}{q(q-1)}\right\}
\end{eqnarray*}
We want to find formulas for the solution 
\[
q^{+}=q(\mathbf{b}^{+}) 
\]
where 
\[
\mathbf{b}^{+}=(b_1,b_2,b_3+1,b_4) 
\]
Let 
\begin{eqnarray*}
h(t) &=&t(t-1)H(p(t),q(t),t,\mathbf{b})+(b_1b_3+b_1b_4+b_3b_4)t \\
&&-\frac 12(b_1b_2+b_1b_3+b_1b_4+b_2b_3+b_2b_4+b_3b_4)
\end{eqnarray*}
and 
\[
h^{+}(t)=h(t)-pq(q-1)+(b_1+b_4)q-\frac 12(b_1+b_2+b_4) 
\]
Following the notation of Okamoto [3], p.354 we have 
\[
A^{+}=(\frac{dh^{+}}{dt}+(b_3+1)^2)(\frac{dh^{+}}{dt}+b_4^2) 
\]
\[
B^{+}=t(t-1)\frac{d^2h^{+}}{dt^2}+(b_1+b_2+b_3+b_4)\frac{dh^{+}}{dt}- 
\]
\[
(b_1b_2(b_3+1)+b_1b_2b_4+b_1(b_3+1)b_4+b_2(b_3+1)b_4) 
\]
\[
C^{+}=2(t\frac{dh^{+}}{dt}-h^{+})-\left\{
b_1b_2+b_1(b_3+1)+b_1b_4+b_2(b_3+1)+b_2b_4+(b_3+1)b_4\right\} 
\]
\[
\]
Then by (2.5) in Okamoto we have the B\"{a}cklund transformation 
\[
q^{+}=\frac 1{2A^{+}}\left\{ (b_3+1+b_4)B^{+}+(\frac{dh^{+}}{dt}%
-(b_3+1)b_4)C^{+}\right\} =\frac UV 
\]
where 
\[
U=t\left\{ p(q-1)(q-t)-(b_3+1)(t-1)-(b_1+b_4)(q-1)\right\} \cdot 
\]
\[
\left\{
pq(q-1)(q-t)-(b_1+b_4)q(q-1)-(b_3+1)q(t-1)+b_1t(q-1)+b_2(q-t)\right\} 
\]
and 
\[
V=(q-t)\left\{ 
\begin{array}{c}
p^2q(q-1)(q-t)^2-p(q-t)(2(b_1+b_4)q^2-(b_1+b_2+2b_4+2b_1t)q+ \\ 
+(b_1+b_2)t)+(b_1^2-(b_3+1)^2)t^2+(b_1b_2+b_1b_4+b_2b_4+(b_3+1)^2)t \\ 
+(b_1+b_4)^2q^2-(b_1+b_4)(2b_1t+b_2+b_4)
\end{array}
\right\} 
\]
The $\tau -$function is defined by (up to a multiplicative constant) 
\[
H(t)=\frac d{dt}\log (\tau (t)) 
\]
After doing the +-construction $n$ times we obtain 
\[
\mathbf{b}_n=(b_1,b_2,b_3+n,b_4) 
\]
\[
q_n 
\]
\[
H_n 
\]
\[
\tau _n 
\]
We have the Toda equation 
\[
\frac d{dt}(t(t-1)\frac d{dt}\log (\tau _n))+(b_1+b_3+n)(b_3+b_4+n)=c(n)%
\frac{\tau _{n+1}\tau _{n-1}}{\tau _n^2} 
\]
where $c(n)$ is a constant which can be chosen to be $1.$

Let $\sigma _n$ be the $n$-th $\tau $-function obtained by replacing $%
\mathbf{b}=(b_1,b_2,b_3,b_4)$ by $\widetilde{\mathbf{b}}%
=(b_1,b_2,b_3,b_4+1). $ Then by (4.16) in Okamoto we have 
\[
q_n=t+\frac{t(t-1)}{b_3+n-b_4}\left\{ \frac d{dt}\log (\sigma _n)-\frac
d{dt}\log (\tau _{n+1})\right\} 
\]

Now we choose as seed solution the rational function (see [2]) 
\[
q=t+\frac{t(t-1)}r\frac{z^{\prime }}z 
\]
where 
\[
z=_2F_1(r,-m,s;t) 
\]
$\mu $ positive integer, which satisfies $P_{VI}(\alpha ,\beta ,\gamma
,\delta )$ with 
\[
\alpha =\frac{r^2}2 
\]
\[
\beta =-\frac{(m+s)^2}2 
\]
\[
\gamma =\frac{(r-s+1)^2}2 
\]
\[
\delta =\frac{1-m^2}2 
\]
This corresponds to 
\[
b_1=\frac{m+r+1}2 
\]
\[
b_2=\frac{m-r+2s-1}2 
\]
\[
b_3=\frac{r+m-1}2 
\]
\[
b_4=\frac{m-r-1}2 
\]
Observe that 
\[
b_1=b_3+1 
\]
We will use that $q$ satisfies the Riccati equation (see [2] ) 
\[
q^{\prime }=\frac 1{t(t-1)}\left\{
(b_4-b_3)q^2+(2b_1t+b_2-b_4-1)q-(b_1+b_2)t\right\} 
\]
Substituting $b_3=b_1-1$ in the formidable expression for $q^{+}$ , it
collapses to 
\[
q^{+}=\frac{t(q-1)\left\{
pq(q-1)(q-t)-(b_1+b_4)q^2+(2b_1+b_4+b_2)q-(b_1+b_2)t\right\} }{(q-t)\left\{
pq(q-1)(q-t)-(b_1+b_4)q^2+(2b_1t+b_2+b_4)q-(b_1+b_2)t\right\} } 
\]
By the Hamiltonian equations we get 
\[
p=\frac{t(t-1)q^{\prime }}{2q(q-1)(q-t)}+\frac 12\left\{ \frac{b_1+b_2}q+%
\frac{b_1-b_2}{q-1}+\frac{b_3+b_4}{q-t}\right\} =\frac{b_1+b_4}{q-t} 
\]
after using the Riccati equation for $q$ . Substituting this into $q^{+}$ we
obtain 
\[
q_1=q^{+}=\frac{(b_1+b_2)t(q-1)}{(2b_1t+b_2-b_1)q-(b_1+b_2)t} 
\]
and 
\[
H_1=H^{+}=H(t)-\frac{pq(q-1)-(b_1+b_4)(q-t)}{t(t-1)}= 
\]
\[
-\frac{(b_1+b_2)(b_1+b_4)}t-\frac{(b_1-b_2)(b_1+b_4)}{t-1} 
\]
which gives 
\[
\tau _1=\exp (\int H_1dt)=\dfrac
1{t^{(b_1+b_2)(b_1+b_4)}(t-1)^{(b_1-b_2)(b_1+b_4)}} 
\]
To start the induction we need also to find $\tau _2.$Since we know $\sigma
_1$(replace $b_4$ by $b_4+1$ in $\tau _1$), namely 
\[
\sigma _1=\frac 1{t^{(b_1+b_2)(b_1+b_4+1)}(t-1)^{(b_1-b_2)(b_1+b_4+1)}} 
\]
we can use the formula for $q_1$%
\[
q_1=t+\frac{t(t-1)}{b_3-b_4+1}\left\{ \frac d{dt}\log (\sigma _1)-\frac
d{dt}\log (\tau _2)\right\} 
\]
i.e. 
\begin{eqnarray*}
\frac d{dt}\log (\tau _2) &=& \\
&&-\frac{(b_1+b_2)(b_1+b_4+1)}t-\frac{(b_1-b_2)(b_1+b_4+1)}{t-1}-\frac{%
(b_1-b_4)(q_1-t)}{t(t-1)}
\end{eqnarray*}
Inspired by numerous experiments we put 
\[
\tau _n=\frac{T_n}{t^{(b_1+b_4)(b_1+b_2+n-1)}(t-1)^{(b_1+b_4)(b_1-b_2+n-1)}} 
\]
so $T_1=1.$ We will later show that the $T_n$ are polynomials. It follows
that 
\[
\frac{T_2^{\prime }}{T_2}=\frac{b_4-b_2}t+\frac{b_4+b_2}{t-1}-\frac{%
(b_1-b_4)(q_1-t)}{t(t-1)} 
\]
In order to compute $T_2$ we have to use the explicit formula for 
\[
q=t+\frac{t(t-1)}r\frac{z^{\prime }}z 
\]
where 
\[
z=_2F_1(r,-m,s,t)=\sum_{j=0}^m(-1)^j\frac{\binom{r+j-1}j\binom mj}{\binom{%
s+j-1}j}t^j 
\]
\textbf{Definition:} 
\[
W(r,m,s)=\sum_{j=0}^m(-1)^j\binom{r+m+1-j}{m-j}\binom{s+m}jt^{m-j} 
\]
\textbf{Lemma 1:} We have the identity 
\[
(s-r-1)z^{\prime }+(r+m+1)(tz^{\prime }+rz)=(-1)^m\frac{(s-1)!m!r(r+1)}{%
(s+m-1)!}W(r,m,s) 
\]
\textbf{Proof:} This is just an identity between binomial coefficients that
is easily verified. 
\[
\]
Expressing everything in the parameters $r,m,s$ and the function $z$ we have 
\[
\frac{T_2^{\prime }}{T_2}=\frac{b_4-b_2}t+\frac{b_4+b_2}{t-1}-\frac{%
(b_1-b_4)(q_1-t)}{t(t-1)}= 
\]
\[
-\frac st+\frac{m-r+s-1}{t-1}-\frac{r+1}{t(t-1)}\frac{%
r(m+s)z-(r+m+1)t((t-1)z^{\prime }+rz)}{(s-r-1)z^{\prime }+(r+m+1)(tz^{\prime
}+rz)}= 
\]
\[
\frac{\left\{ (s-r-1)z^{\prime }+(r+m+1)(z^{\prime }+rz)\right\} ^{\prime }}{%
(s-r+1)z^{\prime }+(r+m+1)(z^{\prime }+rz)} 
\]
\[
\]
after some computations using the hypergeometric equation 
\[
t(1-t)z^{\prime \prime }+(s-(r-m+1)t)z^{\prime }+mrz=0 
\]
Hence we have shown 
\[
\]
\textbf{Proposition 1.. }We have (the constant is of no importance) 
\[
T_2=W(r,m,s) 
\]
Recall that $\sigma _m$ is obtained from $\tau _m$ by the change 
\[
b_4\longrightarrow b_4+1 
\]
In the new parameters this corresponds to 
\[
r\longrightarrow r-1 
\]
\[
m\longrightarrow m+1 
\]
\[
s\longrightarrow s-1 
\]
We define $S_n$ by 
\[
\sigma _n=\frac{S_n}{t^{(b_1+b_4)(b_1+b_2+n-1)}(t-1)^{(b_1+b_4)(b_1-b_2+n-1)}%
} 
\]
Then we have 
\[
S_1=1 
\]
and 
\[
S_2=W(r-1,m+1,s-1) 
\]
The Toda equation for $\tau _n$ and the corresponding one for $\sigma _n$
imply that for $n\geq 2$ we have 
\[
T_{n+1}T_{n-1}=(t^2-t)(T_nT_n^{\prime \prime }-T_n^{\prime \text{ }%
2})+(2t-1)T_nT_n^{\prime }+(n-1)(n+r)T_n^2 
\]
\[
S_{n+1}S_{n-1}=(t^2-t)(S_nS_n^{\prime \prime }-S_n^{\prime \text{ }%
2})+(2t-1)S_nS_n^{\prime }+(n-1)(n+r-1)S_n^2 
\]
Thus we obtain our main result 
\[
\]
\textbf{Theorem:} For $n$,$m$ positive integers we have that 
\[
q_n=t+\frac{t(t-1)}{n+r}\left\{ \frac{S_n^{\prime }}{S_n}-\frac{%
T_{n+1}^{\prime }}{T_{n+1}}-\frac{n+s-1}t-\frac{n+r-m-s}{t-1}\right\} 
\]
satisfies 
\[
P_{VI}(\frac{(n+r)^2}2,-\frac{(m+s)^2}2,\frac{(r-s+1)^2}2,\frac{1-(n+m)^2}2) 
\]
The $T_n$ and $S_n$ are polynomials. 
\[
\]
\textbf{Proof:} To show that the $T_n$ and $S_n$ are polynomials we refer to
the paper [1] where the more general difference equation 
\[
P_{n+1}P_{n-1}=f(x)(P_nP_n^{\prime \prime }-P_n^{\prime \text{ }%
2})+g(x)P_nP_n^{\prime }+h_n(x)P_n^2 
\]
is considered. If 
\[
h_n(x)=n^2+an+b+p(x) 
\]
then 
\[
-2h_{n-1}(x)+h_n(x)=-h_{n-2}(x)+2 
\]
and it follows that the condition 
\[
ff^{\prime \prime }-f^{\prime \text{ }2}+3f^{\prime }g-2fg^{\prime
}-2g^2+2f=0 
\]
is sufficient for the $P_n$ to be polynomials. One checks that $f(x)=x^2-x$
and $g(x)=2x-1$ satisfy this relation.

We have to show that two consecutive $T_n$ are relatively prime. If not then 
$T_n$ and $T_{n-1}$ have a common zero, say $t_0.$ It follows from the
difference equation that then also $T_n^{\prime }(t_0)=0$, i.e. $t_0$ is a
double root of $T_n$ .

Assume first that $t_0\neq 0,1.$ We have 
\[
h(t)=t(t-1)H(t)+\sigma ^{\prime }(\mathbf{b})t-\frac 12\sigma (\mathbf{b}) 
\]
where $\sigma (\mathbf{b})$ is the second symmetric function of $b1,b2,b3,b4$
and $\sigma ^{\prime }(\mathbf{b})$ is the same of $b1,b3,b4.$ Assume that 
\[
h(t)=\frac c{(t-t_0)^k}+\text{lower order terms} 
\]
But $h(t)$ satifies the differential equation 
\[
t^2(t-1)^2h^{\prime }h^{\prime \prime \text{ }2}+\left\{ (2h-(2t-1)h^{\prime
})h^{\prime }+b_1b_2b_3b_4\right\} ^2=\prod_{j=1}^4(h^{\prime }+b_j^2) 
\]
which gives after looking at the highest order terms 
\[
k=1 
\]
\[
c=t_0(t_0-1) 
\]
It follows that 
\[
H(t)=\frac 1{t-t_0}+\text{lower order terms} 
\]
and that after integration that $t_0$ is a simple zero of $\tau (t).$
Contradiction.

To treat the case $t=0$ and $t=1$ we consider the generic case, i.e.we
consider $r$ and $s$ as indeterminates. One sees that $T_n(0)$ and $T_n(1)$
are nonzero polynomials of $r$ and $s$. 
\begin{eqnarray*}
&& \\
&&
\end{eqnarray*}

We can also find a determinantal formula for the $T_n$. Define 
\[
\widetilde{\tau }_n=(t(t-1)^{\tfrac{n(n+r+1)}2}T_{n+1} 
\]
Then we have 
\[
\widetilde{\tau }_0=1 
\]
\[
\widetilde{\tau }_1=W(r,m,s)\cdot (t(t-1))^{\tfrac{r+2}2} 
\]
and 
\[
\delta ^2\log (\widetilde{\tau }_n)=\frac{\widetilde{\tau }_{n+1}\widetilde{%
\tau }_{n-1}}{\widetilde{\tau }_n^{\text{ }2}} 
\]
where 
\[
\delta =t(t-1)\frac d{dt} 
\]
Darboux's formula gives 
\[
\widetilde{\tau }_n=\left| 
\begin{array}{rrrr}
\widetilde{\tau }_1 & \delta \widetilde{\tau }_1 & \cdot \cdot \cdot & 
\delta ^{n-1}\widetilde{\tau }_1 \\ 
\delta \widetilde{\tau }_1 & \delta ^2\widetilde{\tau }_1 & \cdot \cdot \cdot
& \delta ^n\widetilde{\tau }_1 \\ 
\cdot \cdot \cdot & \cdot \cdot \cdot & \cdot \cdot \cdot & \cdot \cdot \cdot
\\ 
\delta ^{n-1}\widetilde{\tau }_1 & \delta ^n\widetilde{\tau }_1 & \cdot
\cdot \cdot & \delta ^{2n-2}\widetilde{\tau }_1
\end{array}
\right| 
\]
This formula should possibly be useful in proving the following 
\[
\]
\textbf{Conjecture 1.}Given 
\[
T_1=1 
\]
\[
T_2=W(r,m,s) 
\]
and 
\[
T_{n+1}T_{n-1}=T_n^2+\frac 1{(n-1)(n+r)}\left\{ (t^2-t)(T_nT_n^{\prime
\prime }-T_n^{\prime \text{ }2})+(2t-1)T_nT_n^{\prime }\right\} 
\]
Then 
\[
q_n=\frac{m+s}{n+r}\cdot \frac{T_n(r,m+1,s)\cdot T_{n+1}(r-1,m,s-1)}{%
T_{n+1}(r,m,s)\cdot T_n(r-1,m+1,s-1)} 
\]
solves 
\[
P_{VI}\left( \frac{(n+r)^2}2,-\frac{(m+s)^2}2,\frac{(r-s+1)^2}2,\frac{%
1-(n+m)^2}2\right) 
\]
\[
\]
\[
\]
The discriminant of $T_n$ defined above factors nicely as a polynomial in $r$
and $s$ .We make the 
\[
\]
\textbf{Conjecture 2.} Define 
\[
h(k,j)=kj^2-\frac{j^3+2j}3 
\]
Then 
\[
discrim(T_n(r,m,s))=const\cdot 
\]
\[
\prod_{j=1}^{m-1}\left\{ (r+n+m-j)(s+j)(s-r+m-1-j)\right\} ^{h(n-1,j)} 
\]
\[
\prod_{j=1}^{m-1}\left\{ (r+1+j)(s+n+m-1-j)(s-r-n+j)\right\} ^{h(m,j)} 
\]
\[
\prod_{j=m}^{n-1}\left\{ (r+1+j)(s+m+m-1-j)(s-r-n+j)\right\} ^{h(j,m)} 
\]
The degree of the discriminant is 
\[
3\binom{m(n-1)}2 
\]
\[
\]

\textbf{Example 1. }We consider the special case when $s=r+2.$ Then one
finds that 
\[
T_n(t)=(t-1)^{m(n-2)}\check{T}_n(t) 
\]
and 
\[
S_n(t)=(t-1)^{(m+1)(n-2)}\check{S}_n(t) 
\]
where 
\[
\deg (\check{T}_n)=m 
\]
and 
\[
\deg (\check{S}_n)=m+1 
\]
One can find explicit formulas for $\check{T}_n$ and $\check{S}_n.$ Define 
\[
V(a,m,b,n)=\sum_{j=0}^b(-1)^{j+1}\binom{n+m+a}{b-j}\binom{a+j}jt^j 
\]
Then we have the following result: 
\[
\]
\textbf{Proposition 2. }We have 
\[
q_n=\frac{V(r+1,m,m+1,n)V(r,m,m,n)}{V(r+1,m,m,n)V(r,m,m+1,n)} 
\]
solves 
\[
P_{VI}(\frac{(n+r)^2}2,-\frac{(m+r+2)^2}2,-\frac 12,\frac{1-(n+m)^2}2) 
\]
\[
\]
\[
\]
\textbf{2.Some numbertheoretic conjectures.}

Experiments suggest that the $T_n$ and $S_n$ contain constant factors
depending on $r$ . 
\[
\]
\textbf{Conjecture 3.} Given 
\[
T_1=1,\text{ }T_2=m!W(r,m,s) 
\]
\[
S_1=1,\text{ }S_2=(m+1)!W(r-1,m+1,s-1) 
\]
Then define $T_n$ and $S_n$ for $n\geq 3$ by 
\[
\frac d{dt}(t(t-1)\frac d{dt}\log (T_n))=(n-1)(n+r)\left\{ \frac{%
T_{n+1}T_{n-1}}{T_n^2}-1\right\} 
\]
\[
\frac d{dt}(t(t-1)\frac d{dt}\log (S_n))=(n-1)(n+r-1)\left\{ \frac{%
S_{n+1}S_{n-1}}{S_n^2}-1\right\} 
\]
for $n\geq 2.$ Then $T_n$ and $S_n$ are in $\mathbf{Z}[r,s,t].$%
\[
\]
\[
\]

In special cases one gets some remarkable difference equations. We give some
examples. 
\[
\]
\textbf{Example 2.} Define $c(n)$ by 
\[
9(n+1)(n+3)\text{ if }n\equiv 1\text{ }mod\text{ 3} 
\]
\[
\frac{(n+1)(n+3)}3\text{ otherwise} 
\]
Let further 
\[
T_1=1,\text{ }T_2=\frac{W(3,2,1)}3=5t^2-5t+1 
\]
\[
S_1=1,\text{ }S_2=W(2,3,0)=(2t-1)(10t^2-10t+1) 
\]
and 
\[
c(n)T_{n+1}T_{n-1}=(t^2-t)(T_nT_n^{\prime \prime }-T_n^{\prime \text{ }%
2})+(2t-1)T_nT_n^{\prime }+(n-1)(n+3)T_n^2 
\]
\[
(n+2)^2S_{n+1}S_{n-1}=(t^2-t)(S_nS_n^{\prime \prime }-S_n^{\prime \text{ }%
2})+(2t-1)S_nS_n^{\prime }+(n-1)(n+2)S_n^2 
\]
\[
\]
for $n\geq 2.$ Then we conjecure that all $T_n$ and $S_n$ have integer
coefficients and $c(n)$ (and $(n+2)^2)$ is best possible, i.e. the
coefficients in the polynomials have no common factor other than one. 
\[
\]
\textbf{Example 3.}Define $c(n)$ by 
\[
\frac{(n+2)(n+4)}8\text{ if }n\text{ is even} 
\]
\[
4(n+2)(n+4)\text{ if }n\equiv 3\text{ mod }4 
\]
\[
16(n+2)(n+4)\text{ if }n\equiv 1\text{ mod 4} 
\]
Let further 
\[
T_1=1 
\]
\[
T_2=\frac{W(4,3,1)}4=(2t-1)(7t^2-7t+1) 
\]
and 
\[
c(n)T_{n+1}T_{n-1}=(t^2-t)(T_nT_n^{\prime \prime }-T_n^{\prime \text{ }%
2})+(2t-1)T_nT_n^{\prime }+(n-1)(n+4)T_n^2 
\]
for $n\geq 2.$ Then we conjecture that $T_n$ has integer coefficients and $%
c(n)$ is best possible. 
\[
\]
\textbf{Example 4.} Define $c(n)$ by 
\[
25(n+3)(n+5)\text{ if }n\equiv 1\text{ mod }5 
\]
\[
\frac{(n+3)(n+5)}5\text{ if }n\equiv 0,2\text{ mod }5 
\]
\[
(n+3)(n+5)\text{ otherwise} 
\]
Let further 
\[
T_1=1 
\]
\[
T_2=\frac{W(5,4,1)}5=42t^4-84t^3+56t^2-14t+1 
\]
and 
\[
c(n)T_{n+1}T_{n-1}=(t^2-t)(T_nT_n^{\prime \prime }-T_n^{\prime \text{ }%
2})+(2t-1)T_nT_n^{\prime }+(n-1)(n+5)T_n^2 
\]
We conjecture that $T_n$ has integer coefficients and $c(n)$ is best
possible. 
\[
\]
Based on these examples we make the 
\[
\]
\textbf{Conjecture 4.} Let $p$ be a prime $\geq 3.$ Define $c(n)$ by 
\[
p^2(p+n)(p+n-2)\text{ if }n\equiv 1\text{ mod }p 
\]
\[
\frac{(p+n)(p+n-2)}p\text{ if }n\equiv 0,2\text{ mod }p 
\]
\[
(p+n)(p+n-2)\text{ otherwise} 
\]
Then 
\[
c(n)T_{n+1}T_{n-1}=(t^2-t)(T_nT_n^{\prime \prime }-T_n^{\prime \text{ }%
2})+(2t-1)T_nT_n^{\prime }+(n-1)(n+p)T_n^2 
\]
for $n\geq 2$ has polynomial $T_n$ with integer coefficients and $c(n)$ is
best possible. 
\[
\]
We have checked the conjecture for 
\[
p=3,5,7,11 
\]
and for $n$ up to $20$ (at least). 
\[
\]
\textbf{Final remark.}After this paper was finished the author found the
polynomials $T_m$ in [4] which up to a factor and some notation agree with
our $T_n$ . There is even a conjectured explicit formula for them (
Conjecture 3.5 ). 
\[
\]
\subsection*{References:} 
\begin{enumerate}
\item G.Almkvist, Polynomial solutions to difference equations connected to
Painlev\'{e} II-VI, CA/0208244.

\item G.Almkvist, Some rational solutions to Painlev\'{e} VI.

\item K.Okamoto, Studies on the Painlev\'{e} equations I, Sixth equation P$%
_{VI}, $Ann. Mat. Pura 146 (1987), 337-381.

\item M.Noumi, S.Okada, K.Okamoto, H.Umemura, Special polynomials associated
with the Painlev\'{e}
equations II in: Integrable systems and algebraic geometry, Kobe-Kyoto 1997,
World Science Publishing 1998.
\end{enumerate}
\end{document}